\documentclass{proc-l}

\usepackage{epsfig}

\newtheorem{theorem}{Theorem}[section]
\newtheorem{lemma}[theorem]{Lemma}

\theoremstyle{definition}
\newtheorem{definition}[theorem]{Definition}

\theoremstyle{remark}
\newtheorem{remark}[theorem]{Remark}

\numberwithin{equation}{section}

\def\Rn{{\mathbb R}^n}
\def\R{{\mathbb R}}
\def\lqq{\lq\lq}
\def\rqq{\rq\rq} 
\def\abrv{.\hskip 0.1truecm}

\begin{document}
{
\title{Distortion in the spherical metric under quasiconformal mappings }
\author{ Peter A. H\"ast\"o }
\address{ Department of Mathematics, University of Helsinki, P.O. Box 4, 00014, Helsinki, Finland}
\email{peter.hasto@helsinki.fi}
\thanks{Supported in part by The Academy of Finland, Research Contract 12132. I would also like to thank M. Vuorinen for suggesting this topic to me.}


\date{June 6, 2001}
\commby{XXX}

\begin{abstract}
This paper contains bounds for the distortion in the spherical metric, that is to say bounds for the constant of H\"older continuity of mappings $f \colon (\Rn,q) \to (\Rn, q)$ where $q$ denotes the spherical metric. The mappings considered are $K$-quasiconformal ($K\ge 1$) and satisfy some normalizations or restrictions. All bounds are explicit and asymptotically sharp as $K \to 1$.
\end{abstract}

\maketitle

}
 
 %%%%%%%%%%%SECTION 1
 %%%%%%%%%%%SECTION 1
 %%%%%%%%%%%SECTION 1
{  
\section{Introduction}

In this paper we derive explicit asymptotically sharp estimates for the 
constant of H\"older continuity in the spherical metric of quasiconformal (qc) mappings. 
These results are based on reducing the problem
to the Euclidean metric case, where H\"older continuity is well-known. 
The paper \cite{FV} 
provides explicit estimates with some nice properties (cf\abrv Lemma \ref{M1FV}) 
that will be shown to hold for the constants in the spherical metric also.

A similar, but more restricted, result has previously appeared as Lemma 4.1 of \cite{Bo} (cited as Lemma \ref{BonLem} in this paper). 
This result gives an estimate for the H\"older constant with respect to the spherical 
metric of planar quasiconformal mappings, which, however, is not asymptotically sharp
as $K\to 1$. 

This paper is organized as follows: in the next section the main results are stated, 
in Section 3 the standard notation is reviewed and some explicit estimates pertaining
to qc mappings are presented. Section 4 contains various lemmas and the 
proof of the main results.
In Section 5, an explicit bound for the H\"older constant in the spherical metric of 
$K$-quasisymmetric real functions is derived.

 %%%%%%%%%%%SECTION 2
 %%%%%%%%%%%SECTION 2
 %%%%%%%%%%%SECTION 2

\section{Main results}

We first introduce some notation and terminology for the formulation of our main results.
For a more complete account, see the following section.

\begin{definition} \label{HolCDef} 
We define $M_i$ (for $R\ge 1$) to be the smallest constant such that the following
inequalities hold for all $x,y\in B^n(R)$ and for all $K$-quasiconformal 
mappings $f\colon\Rn \to \Rn$, with $n \ge 2$, $f(0)=0$ that satisfy the additional
conditions indicated: 

\begin{itemize}
\item[] $\vert f(x) - f(y) \vert \le M_1(K,n,R) \vert x - y\vert^{\alpha},
\ f(B^n) \subset B^n,$
\item[] $\vert f(x) - f(y) \vert \le M_2(K,n,R) \vert x - y\vert^{\alpha}, 
\ f(1)=1,$
\item[] $q(f(x),f(y)) \le M_3(K,n,R) q(x,y)^{\alpha},
\ f(B^n) = B^n,$
\item[] $q(f(x),f(y)) \le M_4(K,n,R) q(x,y)^{\alpha},
\ f(1)=1,$
\end{itemize}
where $\alpha:=K^{1/(1-n)}$ and $q$ is the spherical metric (defined in(\ref{qmetric})).
\end{definition}

\begin{remark} For brevity, we will write $M_1(K,n) := M_1(K,n,1)$, 
$M_2(K,n) := M_2(K,n,1)$, $M_3(K,n) := M_3(K,n,\infty)$ and 
$M_4(K,n) := M_4(K,n,\infty)$. The two last definitions are to be understood
as replacing the condition \lqq $x,y\in B^n(R)$\rqq\  by the condition 
\lqq $x,y\in \overline{\Rn} $\rqq.
\end{remark}

\begin{remark} We can also define analogous constants for $n=1$, see Section 5. 
\end{remark}

\begin{remark} Note that for $M_3(K,n)$ we require that $f(B^n) = B^n$, whereas 
for $M_1(K,n)$ we only have the normalization $f(B^n) \subset B^n$. 
The reason for this discrepancy will become clear in the proof of Theorem \ref{M_3(2)}. 
\end{remark}

\begin{lemma}\label{M1FV} {\rm (\cite{FV}, pp\abrv 115-7)} The constant $M_1(K,n)$
 has the following properties:
\begin{itemize}
\item[{\rm(1)}] $ M_1(K,n)\to 1$ as $ K\to1$,
\item[{\rm(2)}] $ M_1(K,n)$ is bounded for fixed $n$,
\item[{\rm(3)}] $ M_1(K,n)$ is bounded for fixed $K.$
\end{itemize}
\end{lemma}

Moreover, \cite{FV} provides explicit estimates that 
exhibit the described behavior and shows that 
$M_2(K,n)$ has properties (1) and (3), and provides 
an explicit estimate. In terms of these constants we may 
state Bonfert-Taylor's result as

\begin{lemma}\label{BonLem} {\rm (Lemma 4.1, \cite{Bo})} $ M_4(K,2) \le 
128 \cdot 2^{(1-K)/(2K)}$.
\end{lemma}

Note that this constant does not satisfy the most important of the properties above, 
property (1). However, this result proves property (2) for $n=2$, but not
property (1), which is perhaps the most important.
The main results of this paper 
show that (1), (2) and (3) hold for $M_3(K,n)$ 
and (1) and (2) hold for $M_4(K,n)$. \cite{FV} also shows that $M_2(K,n)$ does
not satisfy (3). Unfortunately, we were unable to settle whether $M_4(K,n)$ satisfies 
(2). Our results are stated in the following two theorems:

\begin{theorem}\label{annThm} Let $\lambda_n$ denote the Gr\"otzsch constant 
{\rm (see \cite{AVV}, (8.38))}, $\eta_{K,n}$ be as in Lemma \ref{etaEst} and 
${\hat M}_2(K,n,R)$ be as in Lemma {\ref{annLem}}. Then
\begin{equation}\label{expM_4Est} M_3(K,n)\le M_4(K,n)\le 
{\hat M}_2(K,n,R)m^{2\alpha}\lambda_n^{2(1-\alpha)}(1+0.13(1-\alpha)),\end{equation}
where
\begin{equation} \label{Rdef}
 R:=\sqrt{{1+32^{1-\beta} \over 1-32^{1-\beta}}} > 1,\ \beta:=1/\alpha:=K^{1/(n-1)},\  m:=\eta_{K,n}(1).
\end{equation}
Moreover, ${\hat M}_2(K,n,R)$, and therefore $M_4(K,n)$, is asymptotically sharp for $K\to 1$ with this
choice of $R$. In particular, for $K\le 2$,
$M_4(K,n)$ satisfies the quite crude estimate:
$$ M_4(K,n) \le  m^4\lambda_n^{4(\beta-1)}e^{0.73\sqrt{\beta-1}}(1+0.13(1-\alpha))\le 
\exp \{106\sqrt{K-1} \}. $$
This means that 
$$ q(f(x),f(y)) \le \exp \{138\sqrt{K-1} \} q(x,y)^{\alpha}$$
for $K\le 2$. If we assume instead that $K\le 1.01$ we have the better estimates
$$ M_4(K,n)\le \exp \{ 7 \sqrt{K-1} \},\ 
q(f(x),f(y)) \le \exp \{7\sqrt{K-1} \} q(x,y)^{\alpha}. $$ 
\end{theorem}

\begin{theorem}\label{M_3(2)} The constant $M_3(K,n)$ satisfies (2) in Lemma {\ref{M1FV}}. 
In particular,
$$ M_3(K,n)\le 2^{1-\alpha/2} M_1(K,n)^{1+2\alpha}(1+0.13(1-\alpha)). $$
\end{theorem}

 %%%%%%%%%%%SECTION 3
 %%%%%%%%%%%SECTION 3
 %%%%%%%%%%%SECTION 3

\section{Definitions and preliminary results}

We denote by $\{e_1(=1), e_2(=i), ..., e_n\}$ the 
standard basis of $\Rn$. The following standard notation will be used:
$$ B^n(x, r):=\{y \in \Rn \mid \vert   x-y \vert   <r \},\quad S^{n-1}(x, r):=\partial B^n(x,r),$$ 
$$B^n(r):=B^n(0,r), \quad S^{n-1}(r):=S^{n-1}(0,r), \quad B^n:= B^n(1), \quad
 S^{n-1}:=S^{n-1}(1).$$

We define the spherical metric $q$ by means of the canonical projection onto the
Riemann sphere, hence, for $x,y\in\Rn$, 
\begin{equation}\label{qmetric} q(x,y):={ \vert   x-y \vert   \over  \sqrt{1+\vert  x\vert^2 
}\sqrt{1+\vert y\vert^2 }},\ q(x,\infty):= {1 \over \sqrt{1+\vert  x\vert^2 }}. \end{equation}

We denote the dimension of the Euclidean space under consideration by $n$. 
Till Section 5, it is assumed that $n\ge 2$. The letter $K$ will denote the 
constant of quasiconformality of the arbitrary quasiconformal mapping $f$ 
(of $\overline{\Rn}$). 
More precisely, this entails the following assumptions:
$K \in [1, \infty)$ and $f\colon \overline{\Rn} \to \overline{\Rn}$ 
is an ACL homeomorphism, satisfying
$$ \sup_{\vert  h \vert  =1} \{ \vert f'(x) h \vert \}^n/K 
\le \vert J_f(x) \vert \le
K \inf_{\vert  h \vert =1} \{ \vert f'(x) h \vert \}^n $$
a.e.\ in $\overline{\Rn}$, where $J_f$ stands for the Jacobian 
(cf\abrv \cite{AVV}, Theorem 9.9).
For ease of reference we define the classes 
$$ QC_K:=\{f \colon \overline{\Rn}\to \overline{\Rn} \mid f\ {\rm is\ {\it K}-qc}, 
f(0)=0, f(1)=1, f(\infty)=\infty\}$$
and refer to mappings of this class as normalized $K$-qc mappings. 
Note that for $n\ge 3$,
$QC_1$ consists of isometries only, that is, reflections in planes containing the 
$e_1$-axis and rotations about the same,
by the constancy of the cross ratio under M\"obius mappings. From the 
definition it also follows that if $f \in QC_K$ then $f^{-1} \in QC_K$.

Let $\eta \colon [0,\infty)\to [0,\infty)$ be a homeomorphism. 
Let $D$ and $D'$ be domains in $\Rn$. We say that a mapping $g:D \to D'$ is 
$\eta$-quasisymmetric if
\begin{equation}\label{qsdef} { \vert g(a)-g(c) \vert \over \vert g(a)-g(b) \vert} \le
\eta\left({ \vert a-c \vert \over \vert a-b \vert}\right)\end{equation}
for all $a,b,c \in D$ with $a \neq b$. 

\begin{lemma}\label{etaEst} 
{\rm (\cite{AVV} Theorems 14.6 and 14.8, \cite{Vu} Theorem 7.47)} 
Let $\lambda_n$ as in Theorem \ref{annThm}.  
For $n\ge 2$ and $K\ge 1$ there exists a bijection $\eta_{K,n} \colon [0,\infty)\to [0,\infty)$ such that every $K$-qc mapping of $\Rn$ is $\eta_{K,n}$-qs, 
and $\eta_{K,n}$ satisfies the following bounds: 
{\rm
\begin{itemize}
\item[(1)] $m=\eta_{K,n}(1)\le \exp\{4K(K+1)\sqrt{K-1}\},$
\item[(2)] $\eta_{K,n}(t)\le m{\lambda_n}^{1-\alpha}t^{\alpha}$, for $0\le t \le 1$,
\item[(3)] $\eta_{K,n}(t)\le m{\lambda_n}^{\beta-1}t^{\beta}$, for $t \ge 1$,
\item[(4)] ${\eta_{K,n}}^{-1}(t)\ge {\lambda_n}^{1-\beta}(t/m)^{\beta}$, 
for $0\le t \le m$,
\item[(5)] ${\eta_{K,n}}^{-1}(t)\ge {\lambda_n}^{\alpha-1}(t/m)^{\alpha}$, for $t \ge m.$
\end{itemize}}
\end{lemma}

For the remainder of this paper 
let us write $\alpha:=K^{1/(1-n)}=:1/\beta$ and $m:=\eta_{K,n}(1)$.
It is well-known that $\lambda_n^{1-\alpha}\le 2^{1-1/K}K$ and that 
$\lambda_n\in [4,2e^{n-1}]$ (see \cite{AVV}, p. 169).

\label{etaEstImp}\begin{lemma} {\rm (\cite{AVV}, 14.36 (4))} For $K \le 1.01$, we have the better estimate
$$ m=\eta_{K,n}(1) \le (K-1)^{1-K}e^{9(K-1)}. $$
\end{lemma}

%%%%%%%%%%% SECTION 4
%%%%%%%%%%% SECTION 4
%%%%%%%%%%% SECTION 4

\section{H\"older continuity in the spherical metric}

The main results, Theorems \ref{annThm} and \ref{M_3(2)}, will be proved in this section. 
To accomplish this, we first temporarily turn to a more general setting. 
Suppose, therefore,  that $g$ is an arbitrary ${\gamma }$-biH\"older 
continuous mapping with constant $M$ (with respect to the Euclidean metric). 
What can we say about the H\"older 
continuity of $g$ with respect to the spherical metric? In this 
general setting we may state the following lemma, which 
turns out to be asymptotically sharp for qc-mappings:

\begin{lemma}\label{HolderLem} Let for $0<\gamma\le 1$ consider an arbitrary
 $\gamma-$biH\"older continuous mapping, $g:D \to D'$, fixing the origin (which entails
 that $D$ and $D'$ contain it). That is to say, there exists a $M\ge 1$
such that
\begin{equation}\label{biHE} \left({{ \vert  x-y  \vert } \over {M}}\right)^{1/{\gamma }} \le  
  \vert  g(x)-g(y) \vert   \le  M  {\vert  x-y  \vert }^{{\gamma }},\end{equation}
for all $x,y \in D$.
Then $g$ is H\"older continuous with respect to the spherical metric 
as well, in particular:
$$ q(g(x),g(y)) \le  M' q(x,y)^{\gamma },$$
where $M' \le M^{1+2\gamma}(1+0.13(1-\gamma))$. 
\end{lemma}

\begin{proof} By (\ref{qmetric}), the second inequality of (\ref{biHE}) 
is equivalent to
$$ q(g(x),g(y)) \le M \sqrt{(1+\vert x\vert^2)^{\gamma} \over 1+ \vert g(x)\vert^2}
\sqrt{(1+ \vert y\vert^2)^{\gamma} \over 1+ {\vert g(y)\vert}^2} q(x,y)^{\gamma}$$
and hence $q(g(x),g(y)) \le  M c q(x,y)^{\gamma }$, where 
$$ c:=\sup_{x,g} {(1+ \vert x \vert^2)^{\gamma }\over 1+\vert g(x)\vert^2} ,$$
and the supremum is taken over $x\in D$ and over all mappings $g$ that satisfy 
the assumptions. Using the first inequality of (\ref{biHE}) we conclude, by setting
$y=0$ in (\ref{biHE}), that $\vert g(x)\vert^{\gamma }M \ge \vert x \vert.$
Then, since $M\ge 1$ and $0 < {\gamma } \le  1$, we have
$$ c \le  M^{2{\gamma }} \sup_{t\ge 0}{( 1+t^{\gamma })^{\gamma } \over 1+t}.$$
To prove that 
\begin{equation}\label{1.5c} {( 1+t^{\gamma })^{\gamma } \over 1+t}\le 1+0.13(1-\gamma)\end{equation}
we take the logarithmic partial derivative with respect to $\gamma$ 
of both sides. Since the left and right hand side
of (\ref{1.5c}) are equal when $\gamma=1$, we only need to prove that the derivative of the 
left hand side is greater than that of the right hand side for $0<\gamma<1$.
That is, we need to show that 
$$ \log (1+s) + {s\log s \over 1+s} \ge -{0.13 \over 1+0.13(1-\gamma)}$$
where $s:=t^{\gamma}$. The left hand side has an infimum greater than
$-0.1144$ and the right hand side the supremum $-0.13/1.13<-0.115$, so we are done.  \end{proof}

\begin{remark} Note that a mapping can be $\gamma$-biH\"older continuous in the
sense of (\ref{biHE}) for $\gamma <1$ only if $D$ and $D'$ are bounded.
\end{remark}

\begin{lemma}\label{qHcBn} $ M_3(K,n,1)\le M_1(K,n)^{1+2\alpha}(1.13-0.13\alpha).$
\end{lemma}

\begin{proof} Follows directly from Lemma \ref{HolderLem} as $f$ is $\alpha$-H\"older 
continuous in $B^n$ with constant $M_1(K,n)$ by definition.  
\end{proof} 

From the Lemmas \ref{M1FV} and \ref{HolderLem} we derive the properties (1)-(3)
for $M_3(K,n,1)$. It does not provide the H\"older continuity throughout 
space. It is possible to remedy this
shortcoming as will be shown, however, in doing so, we will loose property (2).
Therefore we first give an independent proof of property (2).

\begin{proof}[Proof of Theorem \ref{M_3(2)}.]  
By Lemma \ref{HolderLem} this is clear if $x,y\in \overline{B^n}$ or 
$x,y \in \{B^n\}^c$. Let $x\in B^n$ and $y\in \overline{B^n}^c$ and $f\in QC_K$ be given. 
We will prove that there exists a point $w$ in the intersection 
of the line through $x$ and $y$ with $S^{n-1}$ such that 
$q(x,w)+q(w,y)\le \sqrt{2}q(x,y)$. Let $x'$ be the image of $x$ in the canonical 
projection of $\overline{\Rn}$ onto $S^n$ etc. Then, by definition, 
$\vert x'-w'\vert =q(x,w)$ and so on. Let $r$ be the radius of the image of 
the line containing $x$ and $y$ under the projection. Of the two possible 
$w$'s, we chose the one on the shorter arc joining $x'$ and $y'$.
Then we need to show that
$$ r\sqrt{2-2\cos \alpha} + r\sqrt{2-2\cos \beta-\alpha} \le 2r\sqrt{1-\cos \beta}.$$
Let us differentiate the left-hand-side with respect to $\alpha$. We the find easily that 
the left-hand-side has a maximum at $\alpha=\beta/2$. We then have to show that
$\sqrt{2-2\cos\beta/2} \le \sqrt{1-\cos\beta}$. Now this follows directly from the
equation $\cos\beta = 2\cos^2\beta/2 - 1$.

From Lemma \ref{HolderLem} it follows that
$$ q(f(x),f(w)) \le M_3(K,n,1) q(x,w)^{\alpha},
\ q(f(w),f(y)) \le M_3(K,n,1) q(w,y)^{\alpha}.$$
Here we have used the fact that a continuous mapping has the same H\"older
constant for $B^n$ and $\overline{B^n}$.

We combine the previous estimates and conclude that
$$ q(f(x), f(y))\le M_3(K,n,1) (q(x,w)^{\alpha}+q(w,y)^{\alpha}) \le$$
$$ \le M_3(K,n,1)2^{1-\alpha} (q(x,w)+q(w,y))^{\alpha} \le
2^{1-\alpha} M_3(K,n,1)(\sqrt{2}q(x,y))^{\alpha}= $$
$$ = 2^{1-\alpha/2}M_3(K,n,1) q(x,y)^{\alpha}.$$ 
Thus $ M_3(K,n)\le 2^{1-\alpha/2} M_3(K,n,1)$.
Note that this estimate is not asymptotically sharp.
\end{proof}

We now return to deriving an asymptotically sharp estimate and start by
presenting an auxiliary lemma stating an explicit relationship
between the constant and the domain of definition of the mapping:

\begin{lemma}\label{annLem} For $R\ge 1$ 
\begin{equation}\label{expAnnEst} M_2(K,n,R) \le {\hat M}_2(K,n,R) :=M_1(K,n)m {\lambda_n}^{\beta-1}R^{\beta-\alpha}.\end{equation}
\end{lemma}

\begin{proof} Let $f\in QC_K$ and $x,y\in B^n(R)$. Define an auxiliary mapping $g$ by:
$$ g(x):={\eta_{K,n}(R)}^{-1}f(Rx).$$
Since $g(B^n)\subset B^n$ it follows from the definition of $M_1(K,n)$ that
$$ M_1(K,n)\vert  x/R- y/R\vert^{\alpha} \ge \vert  g(x/R)-g(y/R)\vert =
{\eta_{K,n}(R)}^{-1}\vert f(x)-f(y)\vert , $$
which is equivalent to $\vert f(x)-f(y)\vert\le M_1(K,n)\eta_{K,n}(R) 
R^{-\alpha}\vert x-y\vert^{\alpha}$
and thus, by Lemma \ref{etaEst}, (\ref{expAnnEst}) holds with the constant 
$M_2(K,n,R)$ indicated.  
\end{proof} 

Now we will derive a new version of Lemma \ref{HolderLem}, which incorporates the 
variability of $M_2(K,n,R)$ in the proof and is thus 
specialized to quasiconformal mappings.

\begin{lemma}\label{HolderLemQC} Let $f\in QC_K$. Then $f \vert_{B^n(R)}$ is
H\"older continuous with respect to the spherical metric with constant 
$M_4(K,n,R)$ satisfying
$$ M_4(K,n,R) \le M_2(K,n,R) m^{2\alpha} \lambda_n^{2(1-\alpha)}(1+ 0.13(1-\alpha)).$$
\end{lemma}

\begin{proof} This is proved as Lemma \ref{HolderLem} except that we will use 
different methods of estimating the upper bounc $c$ of $(1+ \vert x \vert^2)^{\alpha} 
/(1+\vert f(x) \vert^2)$. If $\vert f(x)\vert \ge 1$ we use the estimate
$$ \vert x \vert \le \eta_{K,n}(\vert f(x) \vert) \le m{\lambda_n}^{\beta -1} 
\vert f(x)\vert^{\beta}. $$
In this case we easily see that $c\le (m{\lambda_n}^{\beta -1})^{2\alpha}$. 
For $x \in f^{-1}(B^n)$ we have, by the H\"older continuity of $f^{-1}$,
$\vert x \vert \le m \lambda_n^{1-\alpha} \vert f(x) \vert^{\alpha}$. 
Then we proceed exactly as Lemma \ref{HolderLem} to get 
$c \le (m \lambda_n^{1-\alpha})^{2\alpha}(1+ 0.13(1-\alpha))$. It follows that 
$$ c \le \max \{(m{\lambda_n}^{\beta -1})^{2\alpha}, (m \lambda_n^{1-\alpha})^{2\alpha}
(1+ 0.13(1-\alpha)) \} m^{2\alpha}\lambda_n^{2(1-\alpha)} (1+ 0.13(1-\alpha)).$$  
Recall that in Lemma \ref{HolderLem} the H\"older constant in the spherical 
metric was shown to be less than $M c$ so we are done. 
\end{proof}

We are now in a position to derive easily the main result.

\begin{proof}[Proof of Theorem \ref{annThm}.] 
By Lemmas \ref{annLem} and \ref{HolderLemQC} the estimate holds for $x,y 
\in B^n(R)$. Since the inversion $x \mapsto x \vert x \vert^{-2}$ is a q-isometry,
we conclude that the estimate holds also for $x,y \in B^n(1/R)^c$. We now deal with 
the remaining case: $x\in B^n(1/R), y\in B^n(R)^c$. Then 
$$ q(x,y)\ge {R-1/R \over \sqrt{1+R^2} \sqrt{1+1/R^2} }={R^2-1 \over R^2+1}= 32^{1-\beta}\ge 1/M_4(K,n)^{\beta}.$$
The last inequality holds, since we may assume (otherwise inequality (\ref{expM_4Est})
is trivial) $M_4(K,n)^{\beta}\ge \lambda_n^{3(\beta-1)} \ge 32^{\beta-1}$ 
as $\lambda_n\ge 4$. From this it follows immediately that 
$$ M_4(K,n)q(x,y)^{\alpha} \ge 1 \ge q(f(x),f(y)) $$
which concludes the proof of (\ref{expM_4Est}).

For asymptotical sharpness we still need to prove that $M_2(K,n,R)\to 1$ when 
$K\to 1$ for $R$ as given in (\ref{Rdef}), i.e. that $R^{\beta-\alpha}\to 1$ at 
the limit. This follows from the following chain of inequalities:

$$ R^{\beta-\alpha}\le \left({1+32^{1-\beta} \over 1- 32^{1-
\beta}}\right)^{\beta-1} \le
 \left({2 \over 1- 32^{1-\beta}}\right)^{\beta-1} \le \exp(0.73\sqrt{\beta - 1}).$$
Thus we see that the bound of $M_4(K,n)$ is asymptotically sharp.  
\end{proof}

\begin{remark} 
Note how using this latter method does not give us a better explicit bound for
$M_3(K,n)$ than for $M_4(K,n)$. This is a consequence of our operating in Euclidean
space, where the sharper restriction of the mappings on the unit circle is
not much different from only fixing unity as \lqq viewed\rqq\  from near infinity.
Also note that already for quite small values
of $K$, the estimate in Theorem \ref{annThm} is worse that that in 
the proof of Theorem \ref{M_3(2)}.
\end{remark}

\begin{remark} 
If $K$ is close to one, we may derive the following bound for $M_1(K,n,R)$ in terms of 
our new constant $M_3(K,n,R)$:
$$ M_1(K,n,R)\le {M_3(K,n,R)^{1+2\alpha}\over (1+R)^{\alpha}-M_3(K,n,R)^2R^{\alpha}}
(1.13-0.13\alpha).$$
It is valid if $(1/R+1)^{\alpha}>M_3(K,n,R)^2$.
In particular,
$$ M_1(K,n)\le { M_3(K,n,1)\over 2^{\alpha}-M_3(K,n,1)^2}(1.13-0.13\alpha) \le 
{ M_3(K,n)\over 2^{\alpha}-M_3(K,n)^2}(1.13-0.13\alpha)  $$
if $M_3(K,n,1)<2^{\alpha/2}$. These claims are proved exactly as 
Lemma \ref{HolderLem} so the proofs are omitted here. We see that
from bounds on the H\"older constants in the spherical metric we can
derive bounds for the constants in the Euclidean metric, as well as the other
way around.
\end{remark}

%%%%%%%%%%% SECTION 5
%%%%%%%%%%% SECTION 5
%%%%%%%%%%% SECTION 5

\section{An additional result} 

In \cite{LV}, a one parameter family of $K$-quasisymmetric (qs) functions is 
defined as those functions $g\colon \R \to \R$ that satisfy 
$$ {1\over K} \le {g(x+t)-g(x) \over g(x)-g(x-t)} \le K $$
for all $x$ and $t>0$. We will now derive an explicit bound for the 
H\"older constant in the spherical metric of this class of functions.
The proof of this is similar to the existence proof 
given by Lehto and Virtanen (\cite{LV}, p. 57), 
except that the following lemma will be used:

\begin{lemma}\label{LehLem}{\rm (\cite{Le})} Let $f$ be
$K$-qs in the sense of \cite{LV} with constant $K$. 
Then $f$ can be extended to an $L(K)$-qc mapping $ f\colon 
H \to H$ where 
$$ L(K) \le \min \{ K^{3/2} , 2K-1 \}$$
and $H$ denotes the upper half-plane.
\end{lemma}

\begin{remark} We can further extend the mapping of the upper half-plane 
to a mapping of $\R^2$ by reflection.
\end{remark}

We will also need the following geometrical lemma, which states quite simply 
that a point far away from the real axis and from infinity is not near
the unit circle in the spherical metric.

\begin{lemma}\label{auxLem} Let $(x_1,x_2)\in \R^2$ be a point such that $x_2<-1/R$ 
and $\vert x\vert \le R$ for $R>1$ ($x$ is in the shaded area in Figure 1A). 
Let $\pi$ denote the inversion in $S^1(e_2, \sqrt{2})$. Then 
$$\vert \pi(x) \vert \le \sqrt{R^3+R-2 \over R^3+R+2}.$$
\end{lemma}

\begin{proof} The point $x$ lies in a region bounded by a line and a circle, as shown in 
Figure 1A. This means that $\pi(x)$ lies in a region bounded by two spheres, as shown 
in Figure 1B. It is immediately clear that $y$ is as far from the origin as any other
point in the region, so that $\vert \pi(x)\vert \le \vert y\vert$. 
We calculate $\vert y\vert$ from
$$ {R^2\over (R+1)^2} = {1\over (R+1)^2} + \vert y\vert^2 - 2 {\vert y\vert\over R+1}
\cos\alpha$$
and 
$$ {2R^2\over (R^2-1)^2} = {(R^2+1)^2\over (R^2-1)^2} + \vert y\vert^2 - 
2 \vert y\vert {(R^2+1)\over (R^2-1)}\cos\alpha.$$
(Note that here $\alpha$ refers to an angle as indicated in the picture, not to 
$K^{1/(1-n)}$.) We subtract the second equation from the first to get
$$ 2\vert y\vert \cos\alpha= 2R {R-1 \over R^2-R+2}.$$
Now using this in the first equation gives us the desired formula 
$\vert y\vert^2 = (R^3+R-2)/(R^3+R+2)$. 
\end{proof} 

\bigskip
%\centerline{\epsfig{figure=fig8.eps,width=4in, height=2in}}
\includegraphics [width=4in, height=2in]{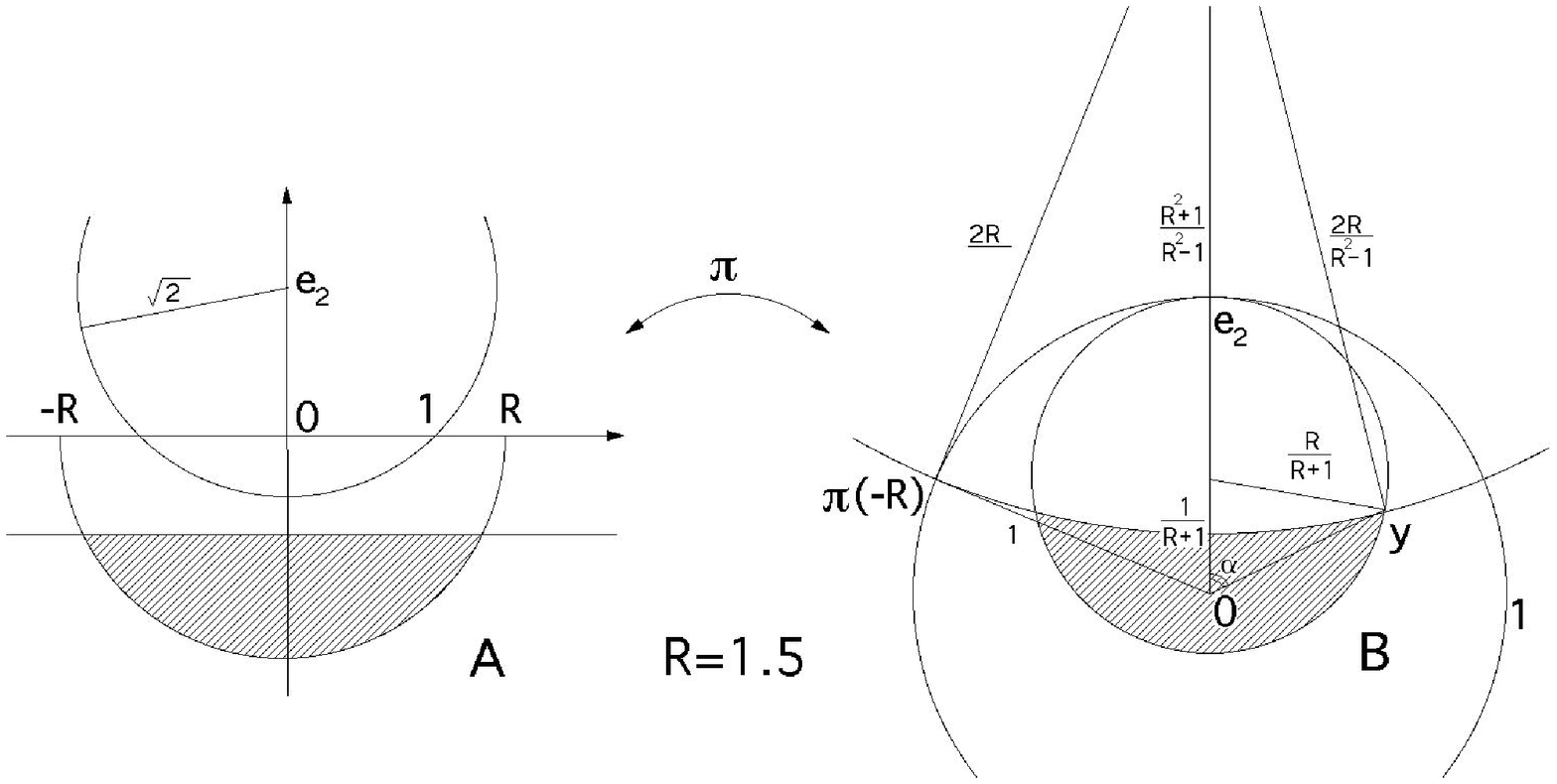}
\\
\centerline{ \bf Figure 1.}
\bigskip

\begin{theorem} Let $f\colon \R \to \R$ be K-qs. Then
$$ q(f(x), f(y)) \le (1/2) 16^{1-1/L} (R^3+R + \sqrt{(R^3+R)^2-4})q(x,y)^{1/L}$$
for all $x,y\in \R$, where $R:= \eta_{L,2}(1)$ and $L\le \min \{K^{3/2}, 2K-1\}$. 
\end{theorem}

\begin{proof} Denote the symmetric $L$-qc extension of $f$ to a mappin from 
$\R^2$ to $\R^2$ by $f$, again (here we use Lemma \ref{LehLem}). We assume without loss of generality that $f(0)=0$ and $f(1)=1$. Since $f$ is $L$-qc, 
$\vert f(-e_2) - f(0) \vert \le \eta_{L,2}(1)$. We may assume without
loss of generality that $f(-e_2)_2 < 0$. Let $x_0$ be a point, such that 
the projection of $f(-e_2)$ in the real axis equals $f(x_0)$. Then
$$\vert f(-e_2)-f(x_0)\vert \ge \eta_{L,2}^{-1}(\vert e_2+x_0\vert) 
\ge \eta_{L,2}^{-1}(1)=1/\eta_{L,2}(1),$$
where the last inequality holds, since $\eta_{L,2}^{-1}$ is increasing and
$\vert e_2+x_0\vert \ge 1$, since $x_0\in \R$. Therefore we may take $R=\eta_{L,2}(1)$
in Lemma \ref{auxLem} to conclude that 
\begin{equation}\label{zero} \vert \pi(f(-e_2))\vert \le \sqrt{R^3+R-2 \over R^3+R+2}\end{equation}
with $\pi$ the inversion in $S^1(e_2, \sqrt{2})$. Let $g:= \pi\circ f\circ \pi$. 
Since $\pi$ is a 1-qc, $g$ is $L$-qc. Since $f$ fixes $\R$, 
$g$ fixes $S^1$. However,  (\ref{zero}) means that $\vert g(0)\vert \le 
\sqrt{R^3+R-2}/\sqrt{R^3+R+2}$, since $\pi(-e_2)=0$. 

Let $\tau$ be a M\"obius transformation with $\tau(S^1)=S^1$ and $\tau(g(0))=0$. Then,
since $g$ is symmetric in $S^1$, $\tau(g(\infty))=\infty$. Now $\tau$ is biLipschitz 
with constant
$$ {1+\vert g(0) \vert\over1-\vert g(0)\vert} \le {\sqrt{R^3+R-2} + \sqrt{R^3+R+2}\over
\sqrt{R^3+R+2} - \sqrt{R^3+R-2}} = {1\over 2} ( R^3+R + \sqrt{(R^3+R)^2-4}).$$
(See e\abrv g\abrv \cite{Vu}, chapter 1, for these elementary M\"obius mapping results.) Since $\tau\circ g$ is $L$-qc and fixes $S^1$, 0 and $\infty$, we have
$ \vert \tau(g(x)) - \tau(g(y))\vert \le 16^{1-1/L} \vert x-y\vert^{1/L}$,
for $x,y\in S^1$ by Theorem 15.5 of \cite{AVV}. It follows from the 
biLipschitz property of $\tau$ that 
$$ \vert g(x) - g(y)\vert \le (1/2) 16^{1-1/L}( R^3+R + \sqrt{(R^3+R)^2-4}) \vert x-y\vert^{1/L}.$$
When we set $x=\pi (x')$ and $y=\pi (y')$, we get
$$ q(f(x'), f(y')) \le (1/2) 16^{1-1/L}( R^3+R + \sqrt{(R^3+R)^2-4}) q(x',y')^{1/L},$$ 
for $x',y' \in \R$, which is what we wanted to show.  
\end{proof}

\begin{remark}\label{sharp} Since
$R \le e^{\pi (L-1/L)}$ by Theorem 10.35 of \cite{AVV}, we readily
get explicit estimates from the above formula. Note that 
$$ (1/2) 16^{1-1/L} (R^3+R + \sqrt{(R^3+R)^2-4}) \to 1$$
as $K\to 1$. 
\end{remark}

\begin{remark}\label{conj} The previous theorem could be extended to higher 
dimensions by using the extension result of \cite{TV} instead of that of 
\cite{Le}, which would allow us to derive estimates similar to those in the begining of 
this paper. This does not, however, give the explicit results that we have
strived for. 
\end{remark}

\bibliographystyle{amsplain}

}\end{document}